# ASYMPTOTICALLY OPTIMAL MULTISTAGE TESTS OF SIMPLE HYPOTHESES

By Jay Bartroff[1]

*University of Southern California*

A family of variable stage size multistage tests of simple hypotheses is described, based on efficient multistage sampling procedures. Using a loss function that is a linear combination of sampling costs and error probabilities, these tests are shown to minimize the integrated risk to second order as the costs per stage and per observation approach zero. A numerical study shows significant improvement over group sequential tests in a binomial testing problem.

**1. Introduction and summary.** Multistage hypothesis tests have practical advantages over fully-sequential tests in many situations since it is often more costly to perform $n$ single experiments than a single experiment of size $n$. The theory of efficient multistage tests has been developed in essentially two directions. The first is general existence and uniqueness results of Schmitz [20], who shows that optimal multistage procedures do exist for a large class of problems and that the optimum has the renewal-type property that at each stage it behaves as if it were starting from scratch given the data so far, and Morgan and Cressie [5, 18], who prove the existence of a multistage competitor of the SPRT. However, these general results do not tell us anything more specific about the optimal tests and certainly not how to apply them without resorting to backward induction-type computer algorithms or artificial truncations. The second direction is truncated (predetermined number of stages) and group sequential (constant stage size) tests, of which many have been developed for clinical trials; see Pocock [19], Wang and Tsiatis [21], Kim and DeMets [12], Eales and Jennison [7, 8], Jennison and Turnbull [11], Barber and Jennison [1] and Lai and Shih [13]. These authors do provide specific tests that successfully address many practical issues arising in clinical trials, but are not concerned with optimality

Received September 2006; revised January 2007.
[1]Supported by NSF Grant DMS-04-03105.

in a general setting, and those that do prove optimality do so under severe restrictions of truncation or constant stage sizes. Lorden [17] presents a three-stage test that has asymptotically the same total sample size as the SPRT and shows that three stages are necessary for any multistage test to have this property.

These previous results do not address a fundamental question in multistage testing: How does one choose the size of the next stage optimally, given the data observed so far and free of oversimplifying restrictions? This paper aims to answer this question by introducing a family of variable stage size multistage tests which can be described by simple, closed-form equations and are asymptotically optimal, without relying on truncations or group sequential restrictions. We focus here on testing simple hypotheses; extension of these ideas to composite hypotheses is discussed in the author's Ph.D. thesis [2].

A common theme in sequential testing is that testing hypotheses can often be reduced to a "power one" test, that is, a test that stops sampling as soon as there is sufficient evidence that the null hypothesis is true but is content to continue sampling forever if it appears that the alternative hypothesis is true. For example, in the fully-sequential setting, Lorden [15, 16] shows that once a substantial number of observations have been taken, asymptotic optimality considerations for testing simple hypotheses can be reduced to considering only power one tests involving the estimated true state of nature versus the opposing hypothesis. Moreover, finding an optimal power one test typically reduces to solving a boundary crossing problem for the relevant test statistic. This suggests the following informal hierarchy:

$$\text{Test of simple hypotheses}$$
$$\textit{reduces to}$$
$$\text{Power one test}$$
$$\textit{reduces to}$$
$$\text{Boundary crossing problem.}$$

In order to derive optimal multistage tests, we consider these three problems in reverse order. In Section 2 we present asymptotically optimal multistage samplers, procedures that sample a random process in stages until it crosses a predetermined boundary. This problem was considered for Brownian motion by Bartroff [3] and we extend those results here to i.i.d., nonnormal data. In Section 3 we use the optimal multistage samplers to design efficient power one tests. In Section 4 we use combinations of these power one tests to design efficient hypothesis tests. Here efficiency is measured by a linear combination of expected sample size, expected number of stages and error probabilities. Our tests are shown to be second order optimal as the costs per stage and per observation approach zero, which corresponds to a large sample size. In marked contrast to constant stage-size group sequential tests, the asymptotically optimal tests and samplers presented here necessarily have stage sizes that decrease roughly as successive iterations of the function $x \mapsto \sqrt{x \log x}$ with probability close to 1, while the average number of stages used is determined by the asymptotics of the ratio of the cost per stage to cost per observation. In Section 5 we propose a finite-sample procedure and present the results of a simulation study comparing it with group sequential tests of hypotheses about the probability of success of Bernoulli trials. The variable stage size tests show substantial improvement over the constant stage size tests.

**2. Multistage samplers.** Consider sampling $X_1, X_2, \ldots$ in stages until $\sum X_i \geq a > 0$ at the end of a stage, and in such a way as to minimize

$$(2.1) \qquad c \cdot EN + d \cdot EM,$$

where $N, M$ are the total sample size and number of stages used. Here $c, d > 0$ represent the costs per observation and per stage, so the sum (2.1) is the average cost incurred in crossing the boundary. On one hand, taking a large number of small stages would make $c \cdot EN$ small but $d \cdot EM$ large; on the other hand, taking a small number of large stages would make $c \cdot EN$ large but $d \cdot EM$ small. Thus, the sampler that minimizes (2.1) can be thought of as the optimal compromise between these two extreme sampling strategies. In this section, after some necessary preliminaries, we define a multistage sampling strategy, show in Theorem 2.1 that it asymptotically minimizes this sampling cost, and show conversely in Theorem 2.2 that any efficient sampler must behave similarly; all theorems are proved in the Appendix. This sampler will be used to construct efficient multistage tests in Sections 3 and 4.

Assume that $X, X_1, X_2, \ldots$ are i.i.d. We say that $X$ is *strongly nonlattice* if the characteristic function $v(t)$ of $X$ satisfies

$$(2.2) \qquad \liminf_{x \to \infty} \left\{ \frac{x^2}{1 - \sup_{\eta \leq t \leq x} |v(t)|} - 2 \log x - 2 \log \log x \right\} > -\infty$$

for some $\eta > 0$. We assume that one of the following three conditions holds:

(2.3)   The distribution of $X$ is strongly nonlattice and $EX^4 < \infty$.

(2.4)   The distribution of $X$ is lattice and $EX^4 < \infty$.

(2.5)   There is an $H > 0$ such that $Ee^{tX} < \infty$ for $|t| < H$.

These conditions are what is needed for the necessary sharp large deviation estimates; see Lemma A.1. We essentially require $X$ to have a finite fourth moment plus to be lattice or strongly nonlattice [(2.3)–(2.4)]. However, if this doesn't hold, then our results are still valid if the moment generating function is finite in a neighborhood of the origin [(2.5)]. Assume that $\mu = EX > 0$. Since the problem is not changed by multiplying the $X_i$ and the boundary $a > 0$ by a positive constant, we assume without loss of generality that $\text{Var}\, X = 1$.

We will describe the stage sizes of a multistage sampler by a sequence of nonnegative integer-valued random variables $N = (N_1, N_2, \ldots)$ such that

$$(2.6) \qquad N_{k+1} \cdot 1\{N_1 + \cdots + N_k = n\} \in \mathcal{E}_n \qquad \text{for all } n \geq 1,$$

where $\mathcal{E}_n$ is the class of all random variables determined by $X_1, \ldots, X_n$. The interpretation of the measurability requirement (2.6) is that by the time $N^k = N_1 + \cdots + N_k$, the end of the first $k$ stages, an observer who knows the values $X_1, \ldots, X_{N^k}$ also knows $N_{k+1}$, the size of the $(k+1)$st stage. We also let $N$ denote the total sample size $N^M$, where $M = \inf\{m \geq 1 : X_1 + \cdots + X_{N^m} \geq a\}$, the total number of stages. A *multistage sampler* is a pair $\delta(x) = (N, M)$, where the argument $x > 0$ is the initial distance to the boundary. When there is no confusion as to which sampler is being used, we will write $S_k = X_1 + \cdots + X_{N^k}$, $S_0 = 0$.

After dividing (2.1) through by $c$, minimizing (2.1) is seen to be equivalent to minimizing

$$(2.7) \qquad EN + h \cdot EM,$$

where $h = d/c$. By Wald's equation,

$$(2.8) \qquad EN = ES_M/\mu = a/\mu + E(S_M - a)/\mu \geq a/\mu,$$

so the sampler that minimizes

$$(2.9) \qquad E(N - a/\mu) + h \cdot EM$$

also minimizes (2.7). Also, using (2.9) instead of (2.7) will lead to a more refined "first-order" asymptotic theory.

The problem of describing the sampler that asymptotically minimizes (2.9) to first-order essentially reduces to considering only certain classes of sequences $\{(a, h)\}$, defined with respect to the *critical functions*

$$(2.10) \quad h_m(x) = x^{(1/2)^m} (\log x)^{1/2 - (1/2)^m} \qquad \text{for } m \geq 1,\ h_0(x) = x.$$

To describe a sampler that asymptotically minimizes (2.9) to first-order, it suffices to consider sequences $\{(a, h)\}$ such that $a \to \infty$. Letting "$\ll$" denote asymptotically of smaller order, it will turn out that good samplers use $m$ stages (with probability approaching 1) if $\{(a, h)\}$ satisfies

$$(2.11) \qquad h_m(a) \ll h \ll h_{m-1}(a)$$

as $a \to \infty$ and use $m$ or $m+1$ stages (with probability approaching 1) if $\{(a,h)\}$ satisfies

(2.12) $$\lim \frac{h}{h_m(a)} \in (0, \infty).$$

A sequence $\{(a,h)\}$ satisfying (2.11) is said to be in the *$m$th critical band*, while one satisfying (2.12) is said to be on the *boundary between critical bands $m$ and $m+1$*. Since it will prove convenient to treat $h$ as a function of $a$, we thus consider (2.9) with $h$ replaced by a function $h(a)$ such that $\{(a,h(a))\}$ is either in the $m$th critical band or on the boundary between critical bands $m$ and $m+1$ (for every sequence of $a$'s approaching $\infty$). That is, let

$$\mathcal{B}_m^o = \{h : (0, \infty) \to (0, \infty) | h_m \ll h \ll h_{m-1}\},$$
$$\mathcal{B}_m^+ = \left\{h : (0, \infty) \to (0, \infty) \Big| \lim_{x \to \infty} h(x)/h_m(x) \in (0, \infty)\right\},$$
$$\mathcal{B}_m = \mathcal{B}_m^o \cup \mathcal{B}_m^+$$

and assume $h \in \mathcal{B}_m$ for some $m \geq 1$. Our notation reflects that, as $a \to \infty$, the average number of stages of an efficient sampler approaches

$$m \quad \text{if } h \in \mathcal{B}_m^o,$$
$$m + \varepsilon \quad \text{if } h \in \mathcal{B}_m^+,$$

where $\varepsilon \in (0, 1)$ is a function of $\lim_{x \to \infty} h(x)/h_m(x)$; Figure 1 summarizes this relationship. We define the *risk* of a sampler $\delta(a) = (N, M)$ to be

(2.13) $$R_h(\delta(a)) = E(N - a/\mu) + h(a)EM.$$

Note that, by (2.8), the definition of risk (2.13) is equivalent to the expectation of a linear combination of the overshoot $S_M - a$ and the number of stages used. Define the Bayes sampler $\delta^* = (N^*, M^*)$ to be one that achieves $R_h^*(a) = \inf_\delta R_h(\delta(a))$.

For $x > 0$ and $z \in \mathbf{R}$, let $t = t(x, z)$ be the unique solution of $(x - \mu t)/\sqrt{t} = z$, that is,

(2.14) $$t(x, z) = x/\mu - \frac{z\sqrt{4x\mu + z^2} - z^2}{2\mu^2}$$

by some simple algebra. Let $z_p$ be the upper $p$-quantile of the standard normal distribution. If the $X_i$ are i.i.d. $N(\mu, 1)$ and $t(x, z_p) = n$ is an integer, then the probability that $X_1 + \cdots + X_n$ exceeds $x$ is $p$. This holds approximately when the $X_i$ are not normal by large deviations and this is why $t$ is useful in parameterizing stage sizes. Let $\Phi$ and $\phi$ denote the standard normal distribution function and density. Let

$$(2.15) \qquad u_m(z) = m + \Phi(z) + \psi^+(z) \cdot \frac{\phi(z)}{1 - \Phi(z)},$$

where $\psi^+(z) = \phi(z) - \Phi(-z)z$ was defined by Chernoff [4]. We extend the domain of $u_m$ to $[-\infty, \infty)$ by adopting the convention $u_m(-\infty) = \lim_{z \to -\infty} u_m(z) = m$. The function $u_m$ appears in the second-order term of the Bayes risk; see Theorem 2.1.

Before defining the asymptotically optimal samplers $\delta^o_{m,h}$ and $\delta^+_{m,z}$, we define an auxiliary sampler $\hat{\delta}_n$ that will be used for the final stages of $\delta^o_{m,h}$ and $\delta^+_{m,z}$. For $n \in \mathbb{N}$, $\hat{\delta}_n$ samples a first stage of size $n$, followed (if necessary) by stages of constant size $\lceil n^{1/2} \rceil$. It is shown in Lemma A.2 in the Appendix that $n = n(a) \to \infty$ can be chosen so that the overshoot of $\hat{\delta}_n$ is not too large but its expected number of stages approaches 1 as $a \to \infty$; for this reason, we refer to $\hat{\delta}_n$ as *bold sampling*.

Finally, we define the samplers $\delta^o_{m,h}$ and $\delta^+_{m,z}$, which are shown to be asymptotically optimal below under different conditions. Namely, the sampler $\delta^o_{m,h}$ will be optimal when $h \in \mathcal{B}^o_m$ and $\delta^+_{m,z}$ will be optimal when $h \in \mathcal{B}^+_m$. These samplers are extensions to nonnormal i.i.d. data of the samplers of Bartroff [3] for Brownian motion. Let $n(x, z) = \lceil t(x, z) \rceil$ and $f(x) =$ $(4/\sqrt{\mu})\sqrt{x\log(x+1)}$. Note that $f^{-1}$ is well defined since $f$ is increasing. The samplers $\delta^o_{m,h}(x)$ are indexed by a positive integer $m$ and a positive function $h$, and the argument $x$ is the initial distance to the boundary. Define $\delta^o_{m,h}$ inductively on $m$ as

$$\delta^o_{1,h}(x) = \hat{\delta}_{n(x,\zeta(x))}(x), \qquad \text{where } \zeta(x) = -(\sqrt{h(x)/x} \wedge \sqrt{(3/2)\log(x+1)}),$$

$$\delta^o_{m+1,h}(x) = \text{1st stage } n(x, \sqrt{(1-2^{-m})\log(x+1)})$$
$$\text{followed (if necessary) by } \delta^o_{m,h\circ f^{-1}}(x - S_1).$$

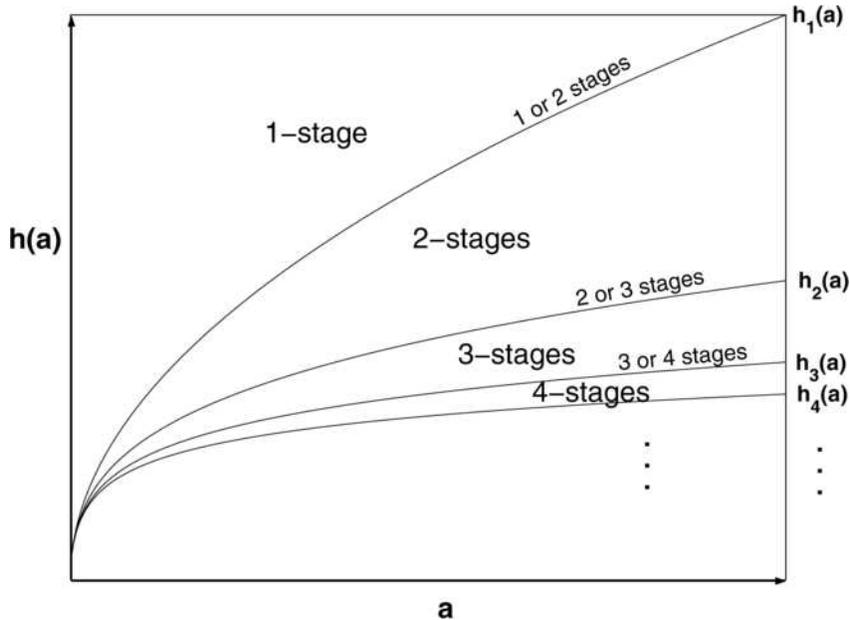

FIG. 1. *The critical functions $h_m$.*

The samplers $\delta^+_{m,z}(x)$, indexed by a positive integer $m$ and a number $z \in \mathbf{R}$, are defined inductively on $m$ as

$$\delta^+_{1,z}(x) = \text{1st stage } n(x,z), \text{ followed (if necessary) by}$$
$$\hat{\delta}_{\nu(x-S_1)}(x-S_1), \text{ where } \nu(y) = n(y, -\sqrt{\log(y+1)}),$$
$$\delta^+_{m+1,z}(x) = \text{1st stage } n(x, \sqrt{(1-2^{-m})\log(x+1)}),$$
$$\text{followed (if necessary) by } \delta^+_{m,z}(x-S_1).$$

THEOREM 2.1. *Assume $h \in \mathcal{B}_m$. Let $z^* \in [-\infty, \infty)$ be the unique solution of*

$$\frac{\phi(z^*)}{1-\Phi(z^*)} = \lim_{x\to\infty} \frac{\kappa_m h_m(x)}{h(x)}, \tag{2.16}$$

*where*

$$\kappa_m = \mu^{-2+(1/2)^m} \prod_{i=1}^{m-1} [(1/2)^{m-1-i} - (1/2)^{m-1}]^{(1/2)^{i+1}}. \tag{2.17}$$

*Then $R^*_h(a) \sim u_m(z^*)h(a)$ as $a \to \infty$. If*

$$\delta = \begin{cases} \delta^o_{m,h}, & \text{if } h \in \mathcal{B}^o_m, \\ \delta^+_{m,z^*}, & \text{if } h \in \mathcal{B}^+_m, \end{cases}$$

*then as $a \to \infty$,*

$$R_h(\delta(a)) \sim u_m(z^*)h(a). \tag{2.18}$$

Theorem 2.2 provides a converse to Theorem 2.1, showing that the type of sampling used by $\delta^o_{m,h}$ and $\delta^+_{m,z}$ is necessary for any efficient procedure. Let $F_y(x) = \sqrt{x\log(y/y^2)}$ and for a function $h$ and $k \in \mathbb{N}$ define

$$F_h^{(k)}(x) = F_y^{(k)}(x)|_{y=h(x)},$$

where the superscript $(k)$ on the right-hand side denotes the $k$th iterate. Bartroff ([3], Lemma 8) showed that $F_h^{(k)}(a)$ is the order of magnitude of how far $\delta_{m,h}^o$ and $\delta_{m,z}^+$ are from the boundary (with probability approaching 1) after the $k$th stage. Theorem 2.2 shows that any sampler that does not follow this "schedule" is necessarily suboptimal.

THEOREM 2.2. *Assume that $h \in \mathcal{B}_m$ and let*

$$\delta = \begin{cases} \delta_{m,h}^o, & \text{if } h \in \mathcal{B}_m^o, \\ \delta_{m,z^*}^+, & \text{if } h \in \mathcal{B}_m^+, \end{cases}$$

*where $z^*$ is as in (2.16). If $\delta' = (N, M)$ is a sampler such that there is a sequence $a_i \to \infty$ with*

(2.19) $\quad P(a_i - S_k \geq (1-\varepsilon)(1/\mu)^{1-2^{-k}} F_h^{(k)}(a_i)) \qquad$ *bounded below 1*

*for some $1 \leq k < m$ and $\varepsilon > 0$, then*

(2.20) $$\frac{R(\delta'(a_i)) - R(\delta(a_i))}{h(a_i)} \to +\infty$$

*as $i \to \infty$. In particular, (2.20) holds if $P(M \geq m) \not\to 1$.*

**3. Power one tests.** Consider the problem of deciding between two densities $f_0$ and $f_1$ by sampling data in stages. Suppose that if $f_0$ is the true density, sampling costs are high and so we want to stop sampling as soon as possible and reject the hypothesis $f_1$. On the other hand, if $f_1$ is the true density, suppose that sampling costs nothing and we are content to observe the data ad infinitum. As an example, suppose a new drug is being marketed under the hypothesis that its side effects are insignificant. Physicians prescribing the drug record and report on the side effects and if they appear unacceptably high $(f_0)$, this must be announced and the drug withdrawn from use. But as long as the hypothesis of insignificant side effects $(f_1)$ remains tenable, no action is required. Although this is an idealized example, power one tests are important theoretical tools because we will use combinations of them to derive optimal hypothesis tests; see Section 1 and the paragraph preceding Section 4.1.

Let $X_1, X_2, \ldots$ be i.i.d. with density either $f_0$ or $f_1$, two distinct densities with respect to some nondegenerate $\sigma$-finite measure. Define a *power one test of $f_0$ versus $f_1$* to be a pair $\delta = (N, M)$ such that $N = (N_1, N_2, \ldots)$ is a sequence of nonnegative integer-valued random variables satisfying the measurability requirement (2.6), with $N_k$, $N^k$ and $M$ defined as in Section 2. Note that a "power one test of $f_0$ versus $f_1$" may only reject $f_1$. If one pays costs per observation and per stage under $f_0$, plus a cost for terminating sampling under $f_1$, then a natural measure of the performance of a power one test of $f_0$ versus $f_1$ is the expected sum of these costs. Hence, we define the *risk* of a power one test $\delta = (N, M)$ of $f_0$ versus $f_1$ to be

$$(3.1) \qquad R_{c,d}(\delta) = cE_0 N + dE_0 M + P_1(N < \infty),$$

where $c, d > 0$. Let $\delta^* = (N^*, M^*)$ be a Bayes test which achieves risk $R^*_{c,d} = \inf_\delta R_{c,d}(\delta)$.

In this section we define a family of power one tests and show in Theorem 3.1 that they minimize the risk to second-order as $c, d \to 0$. Obviously the risk (3.1) depends on the rates at which $c$ and $d$ approach 0, much in the same way that in Section 2 the efficiency of a multistage sampler depended on the asymptotic properties of the function $h$, representing the ratio of the cost per stage to the cost per observation, with respect to the critical functions (2.10). It will turn out that the behavior of efficient hypothesis tests will be determined by an analogous relationship, but with $d/c$ in place of $h$ and a multiple of $\log d^{-1}$ in place of the boundary $a$ in (2.11) and (2.12). That is, it will turn out that efficient hypothesis tests use $m$ stages (with probability approaching 1) if $c, d \to 0$ in such a way that

$$(3.2) \qquad h_m(\log d^{-1}) \ll d/c \ll h_{m-1}(\log d^{-1})$$

and will use $m$ or $m+1$ stages (with probability approaching 1) if

$$(3.3) \qquad \lim_{c,d \to 0} \frac{d/c}{h_m(\log d^{-1})} \in (0, \infty).$$

By analogy with Section 2, we give an essentially complete description of the problem while assuming $c, d \to 0$ at rates satisfying (3.2) or (3.3). To update our notation, let $\mathcal{B}^o_m$ be the set of all sequences $\{(c,d)\}$ such that $1 \geq c, d \to 0$ and satisfying (3.2), let $\mathcal{B}^+_m$ be the set of all such sequences satisfying (3.3), and let $\mathcal{B}_m = \mathcal{B}^o_m \cup \mathcal{B}^+_m$. We prove our main asymptotic results below for sequences $\{(c,d)\} \in \mathcal{B}_m$ for some $m \geq 1$. Note that $\{(c,d)\} \in \mathcal{B}_m$ implies $h_m(\log d^{-1}) = O(d/c)$; hence, a consequence of this assumption is that $d/c \to \infty$. If it were that $d/c$ were bounded below $\infty$, it can be shown that a test with constant stage size and number of stages approaching $\infty$ minimizes the risk (3.1) to second-order. Since our main interest here is variable stage size tests with a small number of stages, we can be sure that the assumption $\{(c,d)\} \in \mathcal{B}_m$ does not exclude any interesting cases.

In this section we use the multistage samplers of Section 2 as power one tests by sampling the log-likelihood process $\log(f_0(X_i)/f_1(X_i))$ until $\sum \log(f_0(X_i)/f_1(X_i))$ exceeds a predetermined boundary. Let

$$(3.4) \qquad \begin{aligned} \sigma^2 &= \mathrm{Var}_0 \log(f_0(X_1)/f_1(X_1)), \\ Y_i &= \sigma^{-1} \log(f_0(X_i)/f_1(X_i)), \end{aligned}$$

so that $E_0 Y_i = \sigma^{-1} I_0 > 0$ and $\text{Var}_0 Y_i = 1$, where $I_0 = E_0 \log(f_0(X_1)/f_1(X_1))$ is the Kullback–Leibler information number. Whenever we use a multistage sampler as a power one test in what follows, we mean with respect to $Y_1, Y_2, \ldots$, which we assume satisfy one of (2.3)–(2.5). Our main result in this section is that the asymptotically optimal multistage samplers derived in Section 2 are second-order optimal as power one tests.

THEOREM 3.1. *Assume that the $Y_i$ satisfy one of (2.3)–(2.5), $\{(c,d)\} \in \mathcal{B}_m$, and let $z^* \in [-\infty, \infty)$ be the unique solution of*

$$(3.5) \qquad \frac{\phi(z^*)}{1-\Phi(z^*)} = \lim_{c,d \to 0} \frac{\kappa_m(\sigma^{-1} I_0) h_m(\sigma^{-1} \log d^{-1})}{d/c},$$

*where $k_m(\mu)$ is as in (2.17). Then $R_{c,d}^* = c I_0^{-1} \log d^{-1} + u_m(z^*) d + o(d)$ as $c, d \to 0$. If $\delta$ is the power one test*

$$\delta = \begin{cases} \delta_{m,d/c}^o(\sigma^{-1} \log d^{-1}), & \text{if } \{(c,d)\} \in \mathcal{B}_m^o, \\ \delta_{m,z^*}^+(\sigma^{-1} \log d^{-1}), & \text{if } \{(c,d)\} \in \mathcal{B}_m^+, \end{cases}$$

*then as $c, d \to 0$,*

$$(3.6) \qquad R_{c,d}(\delta) = c I_0^{-1} \log d^{-1} + u_m(z^*) d + o(d).$$

**4. Tests of simple hypotheses.** In this section we use the optimal power one tests from the previous section to derive optimal multistage tests of two simple hypotheses. Consider the problem of deciding between two distinct densities $f_0$ and $f_1$ by sampling the i.i.d. $X_1, X_2, \ldots$ in stages, while incurring a cost per observation $c$, a cost per stage $d$ and a penalty $w_i$ for incorrectly rejecting $f_i$. Specifically, a *test* of the hypotheses $H_0: f_0$ versus $H_1: f_1$ is a triple $\delta = (N, M, D)$, where $N, M$ are as in Section 3 and $D$ is the "decision" variable taking values in $\{0, 1\}$. The event $\{D = i\}$ means rejection of $H_{1-i}$. Define the *integrated risk* of a test $\delta = (N, M, D)$ with respect to the prior $\pi$ to be

$$r_{c,d}(\delta) = \sum_{i=0}^{1} \pi_i [c E_i N + d E_i M + w_i P_i(D = 1 - i)],$$

where $\pi_i, c, d, w_i > 0$. Let $\delta^* = (N^*, M^*, D^*)$ denote a Bayes test, one that achieves integrated risk $r_{c,d}^* = \inf_\delta r_{c,d}(\delta)$. In this section we define a family of tests and show in Theorems 4.1 and 4.2 that they minimize the integrated risk to second-order as $c, d \to 0$. Moreover, the proofs of these results in the Appendix show that the integrated risk of efficient procedures is dominated by sampling and staging costs; hence, this Bayesian setup can be thought of as a stepping stone to finding tests that are efficient in the frequentist sense as well.

As in Section 3, we assume that $c, d \to 0$ at rates such that $\{(c, d)\} \in \mathcal{B}_m$ for some $m \geq 1$. For $i = 0, 1$, let $\sigma_i^2 = \text{Var}_i \log(f_i(X_1)/f_{1-i}(X_1))$ and $Y_j^{(i)} = \sigma_i^{-1} \log(f_i(X_j)/f_{1-i}(X_j))$ for $j = 1, 2, \ldots$ so that $E_i Y_j^{(i)} = \sigma_i^{-1} I_i$ and $\text{Var}_i Y_j^{(i)} = 1$, where $I_i = E_i \log(f_i(X_1)/f_{1-i}(X_1))$. Whenever we speak of a power one test of $f_i$ versus $f_{1-i}$ (i.e., a test which can only reject $H_{i-1} : f_{1-i}$) below, we will always mean the one defined with respect to $Y_1^{(i)}, Y_2^{(i)}, \ldots$, which we assume satisfy one of (2.3)–(2.5). Let $l_n = \prod_{i=1}^n (f_0(X_i)/f_1(X_i))$ denote the likelihood ratio, and when there is no confusion which $N$ we are considering, we will let $l^k = l_{N^k}$.

To describe the family of optimal tests, we must consider separately two cases of the relationship between $f_0$ and $f_1$. The first case, considered in Section 4.1, is when $I_0 = I_1$ and $\text{Var}_0 X_i = \text{Var}_1 X_i$. This is the "symmetric" case in the sense that the two corresponding power one tests dictate the same initial stage size, and hence, their first stages can be applied simultaneously. This case is of interest because it contains, most notably, the Normal mean problem, $H_0 : \mu = \mu_0$ versus $H_1 : \mu = \mu_1$, about the mean $\mu$ of Normal random variables with known variance, and the symmetric Binomial case, $H_0 : p = 1/2 - \Delta$ versus $H_1 : p = 1/2 + \Delta$, about the probability $p$ of success of a Bernoulli trial. If $I_0 \neq I_1$, the nature of a Bayes test is fundamentally different. In this case, considered in Section 4.2, the ratio of the two initial stages given by the power one tests does not tend to 1, and it is not obvious what the size of the initial stage should be. This gives rise to a necessary "exploratory" first stage, equal to the smaller of the two initial stages dictated by the two corresponding power one tests. The remaining case, where $I_0 = I_1$ and $\text{Var}_0 X_i \neq \text{Var}_1 X_i$, is at present unsolved, but the popular examples contained in the former and the generality of the latter make our analysis sufficient for most purposes.

For simplicity, we present our results here for tests of two simple hypotheses, but these methods and results generalize immediately to tests of $s \geq 2$ simple hypotheses. The asymptotically optimal test for $s > 2$ or for either subcase considered below for $s = 2$ may be loosely described as follows: Sample at the first stage the size of the smallest first stage of the corresponding $s(s-1)$ power one tests, then continue sampling with the power one test of the most likely hypothesis versus the second most likely, according to the results of the first stage.

4.1. *Case I*: $I_0 = I_1$ *and* $\text{Var}_0 X_i = \text{Var}_1 X_i$. Let $(N^{(0)}, M^{(0)})$ be the power one test of $f_0$ versus $f_1$ defined in Theorem 3.1 and let $(N^{(1)}, M^{(1)})$ be the corresponding power one test of $f_1$ versus $f_0$. Under the assumptions $I_0 = I_1$ and $\text{Var}_0 X_i = \text{Var}_1 X_i$, the two procedures $(N^{(0)}, M^{(0)})$ and $(N^{(1)}, M^{(1)})$ dictate the same first stage size. Define the first stage of $\delta = (N, M, D)$ to be this common first stage size, $N_1 = N_1^{(0)} = N_1^{(1)}$. If $l_{N_1} \geq 1$, continue with $(N^{(0)}, M^{(0)})$, stopping the first time $l_{N^k} \geq d^{-1}$ to reject $H_1$, as dictated by $(N^{(0)}, M^{(0)})$, or $l_{N^k} \leq d$ to reject $H_0$. Otherwise, $l_{N_1} < 1$, so switch and continue sampling with $(N^{(1)}, M^{(1)})$ with the same stopping rule. This test is second-order asymptotically optimal, recorded as Theorem 4.1.

THEOREM 4.1. *If the $Y_j^{(i)}$ satisfy one of (2.3)–(2.5), $I_0 = I_1$, $\mathrm{Var}_0 X_i = \mathrm{Var}_1 X_i$, and $\{(c,d)\} \in \mathcal{B}_m$, then*

$$r_{c,d}^* = cI_0^{-1} \log d^{-1} + u_m(z^*)d + o(d), \tag{4.1}$$

$$r_{c,d}(\delta) = cI_0^{-1} \log d^{-1} + u_m(z^*)d + o(d) \tag{4.2}$$

*as $c, d \to 0$, where $z^* \in [-\infty, \infty)$ is the unique solution of*

$$\frac{\phi(z^*)}{1 - \Phi(z^*)} = \lim_{c,d \to 0} \frac{\kappa_m(\sigma_0^{-1} I_0) h_m(\sigma_0^{-1} \log d^{-1})}{d/c}.$$

4.2. *Case II: $I_0 \neq I_1$.* Assume $I_0 < I_1$. For $i = 0, 1$ let $\delta_i(c, d, z)$ denote the power one test of $f_i$ versus $f_{1-i}$ (i.e., the test that can only reject $f_{1-i}$) defined in Theorem 3.1 with generic parameters $c, d, z$. Given $\{(c,d)\} \in \mathcal{B}_m$, define $z_0^*, z_1^*$ to be the unique solutions of the equations

$$\frac{\phi(z_0^*)}{1 - \Phi(z_0^*)} = \lim_{c,d \to 0} \frac{\kappa_m(\sigma_0^{-1} I_0) h_m(\sigma_0^{-1}[1 - I_0/I_1] \log d^{-1})}{d/c}, \tag{4.3}$$

$$\frac{\phi(z_1^*)}{1 - \Phi(z_1^*)} = \lim_{c,d \to 0} \frac{\kappa_m(\sigma_1^{-1} I_1) h_m(\sigma_1^{-1} \log d^{-1})}{d/c}. \tag{4.4}$$

Define $\delta = (N, M, D)$ as follows: Let the first stage of $\delta$ equal the

1st stage of $\delta_1(c, d, z_1^*) = \min_i \{\text{1st stage of } \delta_i(c, d, z_i^*)\}.$

After the first stage,

if $l^1 < 1$, continue sampling with $\delta_1(c, d, z_1^*)$,

if $l^1 \geq 1$, switch and continue sampling with $\delta_0(l^1 c, l^1 d, z_0^*)$,

with the stopping rule

(4.5) stop after the $k$th stage and reject $H_0$ if $l^k \leq d$,

(4.6) stop after the $k$th stage and reject $H_1$ if $l^k \geq d^{-1}$.

Note that $\delta$ stops no later than whichever power one test it chooses after the first stage since $\delta_1(c, d, z_1^*)$ stops when $\sum_1^{N^k} Y_j^{(1)} \geq \sigma_1^{-1} \log d^{-1}$, which is equivalent to (4.5), while $\delta_0(l^1 c, l^1 d, z_0^*)$ stops when $\sum_{N_1+1}^{N^k} Y_j^{(0)} \geq \sigma_0^{-1} \log(l^1 d)^{-1}$, which is equivalent to (4.6). However, $\delta$ may stop before the corresponding power one test because of the stopping rule (4.5)–(4.6). Theorem 4.2 establishes the second-order optimality of $\delta$.

THEOREM 4.2. *If the $Y_j^{(i)}$ satisfy one of (2.3)–(2.5), $I_0 < I_1$, and $\{(c,d)\} \in \mathcal{B}_m$, then*

$$(4.7) \qquad r_{c,d}^* = \sum_{i=0}^{1} \pi_i \{cI_i^{-1} \log d^{-1} + d[1 - i + u_m(z_i^*)]\} + o(d),$$

$$(4.8) \qquad r_{c,d}(\delta) = \sum_{i=0}^{1} \pi_i \{cI_i^{-1} \log d^{-1} + d[1 - i + u_m(z_i^*)]\} + o(d)$$

*as $c, d \to 0$, where $z_i^*$ is given by (4.3) and (4.4).*

**5. A numerical example.** The tests proved asymptotically optimal in Theorems 4.1 and 4.2 are asymptotic not only in the sense that their optimality is proved in the limit as $c, d \to 0$, but also in that they are defined in terms of the rates at which $c, d \to 0$. Thus, in practice, there may be more than one asymptotically optimal procedure for a statistician to choose from. In this section we describe one such procedure and give the results of a numerical experiment comparing it to a sampling with constant stage size.

Given values $0 < c, d < 1$, let

$$(5.1) \quad m_i^* = \inf\{m \geq 1 : \kappa_m(\mu_i) h_m(a_i + 1) - \kappa_{m+1}(\mu_i) h_{m+1}(a_i + 1) \leq d/c\}$$

for $i = 0, 1$, where $\mu_i = I_i/\sigma_i$ and $a_i = \sigma_i^{-1} \log d^{-1}$. Let $\delta$ be the test whose first stage is the smaller of the two first stages of the samplers $\delta_{m_i^*, d/c}^o(\sigma_i^{-1} \log d^{-1})$, and then continues sampling according to

$$(5.2) \qquad \begin{aligned} \delta_{m_0^*, d/c}^o(\sigma_0^{-1} \log d^{-1}) & \qquad \text{if } l^1 > 1, \\ \delta_{m_1^*, d/c}^o(\sigma_1^{-1} \log d^{-1}) & \qquad \text{if } l^1 \leq 1. \end{aligned}$$

The test $\delta$ is asymptotically optimal by Theorem 4.1 when $c, d \to 0$ such that $\{(c,d)\} \in \mathcal{B}_m^o$ since, clearly, $m_i^*$ will equal $m$ for sufficiently small $c, d$.

We consider testing the hypotheses $H_0 : p = 0.4$ versus $H_1 : p = 0.6$ about the probability $p$ of success of i.i.d. Bernoulli trials. To isolate the effects of using variable stage sizes, we compare $\delta$ with the test $\delta_k$ that uses stage sizes of constant size $k$ but with the same stopping rule (4.5)–(4.6), that is, stop when the log-likelihood exceeds $\log d^{-1}$ in absolute value. Table 1 contains the expected sample size, expected number of stages and integrated risk of $\delta$ and $\delta_k$ for various $k$, $c$ and $d$, each of which is computed by 100,000 Monte Carlo replications. For each value of $d/c$, the operating characteristics of $\delta_k$ are given in Table 1 for the following five values of $k$: $k = 1$ (fully-sequential sampling), the (rounded) "average stage size" $EN/EM$ of $\delta$, the size of the first stage of $\delta$, the (rounded) expected sample size $EN$ of $\delta$ and the optimal value $k = k^*$ minimizing $r_{c,d}(\delta_k)$, found by exhaustion. Here $E(\cdot)$ denotes $\sum_{i=0}^{1} \pi_i E_i(\cdot)$. Since both $\delta$ and $\delta_k$ sample until the absolute value of the log-likelihood ratio exceeds $\log d^{-1}$, the cost of the average number of observations required to do this and the cost of the first stage represent "fixed costs," which it is shown in Lemma A.4 in the Appendix that any efficient test must incur. We obtain a more accurate comparison of the efficiency due to variable stage size sampling by considering the *second-order risk* $r'_{c,d} = r_{c,d} - (cEN^{(1)} + d)$, where $N^{(1)}$ is the sample size of $\delta_{k=1}$. The fifth column of Table 1 contains the second-order risk and its percent decreases by $\delta$ in the sixth column. Also included in Table 1 are the asymptotic approximations $m_i^*$ [given by (5.1)] of the optimal expected number of stages and $\widetilde{r}_{c,d} = c \log d^{-1}/I + m_i^* d$ of the Bayes integrated risk.

TABLE 1
*Expected sample size, number of stages, integrated risk and 2nd order risk of $\delta$ and $\delta_k$ for the binomial testing problem $p = 0.4$ vs. $p = 0.6$ with $\log d^{-1} = 10$, $\pi_i = 1/2$, $w_i = 1$*

| Test | $EN$ | $EM$ | $r_{c,d}/d$ | $r'_{c,d}/d$ | $1 - r'_{c,d}(\delta)/r'_{c,d}$ |
|---|---|---|---|---|---|
| | | | $d/c = 5$ | | |
| | | $m_i^* = 7.0$ | $\widetilde{r}_{c,d}/d = 31.7$ | | |
| $\delta$ | 125.1 | 5.9 | 30.9 | 5.33 | – |
| $\delta_1$ | 124.9 | 124.9 | 151.0 | 125.43 | 95.8% |
| $\delta_{21}$ | 141.3 | 6.7 | 35.4 | 9.83 | 45.8% |
| $\delta_{59}$ | 163.3 | 2.8 | 35.9 | 10.33 | 48.3% |
| $\delta_{125}$ | 189.0 | 1.5 | 39.7 | 14.13 | 62.2% |
| $\delta_{k^*=43}$ | 154.2 | 3.6 | 34.4 | 8.83 | 39.6% |
| | | | $d/c = 10$ | | |
| | | $m_i^* = 5.0$ | $\widetilde{r}_{c,d}/d = 17.3$ | | |
| $\delta$ | 130.9 | 4.3 | 17.4 | 4.07 | – |
| $\delta_1$ | 125.1 | 125.1 | 138.5 | 125.17 | 96.7% |
| $\delta_{30}$ | 149.3 | 5.0 | 20.0 | 6.67 | 39.0% |
| $\delta_{70}$ | 171.2 | 2.4 | 20.1 | 6.77 | 39.9% |
| $\delta_{130}$ | 195.4 | 1.5 | 21.0 | 7.67 | 46.9% |
| $\delta_{k^*=49}$ | 157.7 | 3.2 | 19.0 | 5.67 | 28.2% |
| | | | $d/c = 25$ | | |
| | | $m_i^* = 2.0$ | $\widetilde{r}_{c,d}/d = 6.9$ | | |
| $\delta$ | 144.4 | 2.6 | 8.36 | 2.43 | – |
| $\delta_1$ | 124.9 | 124.9 | 131.0 | 125.07 | 98.1% |
| $\delta_{56}$ | 163.8 | 2.9 | 9.92 | 3.99 | 39.1% |
| $\delta_{89}$ | 176.8 | 2.0 | 9.06 | 3.13 | 22.4% |
| $\delta_{144}$ | 201.0 | 1.4 | 9.66 | 3.73 | 34.9% |
| $\delta_{k^*=95}$ | 178.8 | 1.9 | 9.04 | 3.11 | 21.9% |

The results show that $\delta$ has substantially smaller risk and second-order risk than the $\delta_k$. Since $\delta$ and $\delta_k$ use the same stopping rule, this is due to the variable stage sizes of $\delta$, versus the constant stage sizes of $\delta_k$. Even when compared to $\delta_{k^*}$ with the optimal fixed stage size $k^*$ (which requires fitting an additional parameter), $\delta$ has roughly 40%, 30% and 20% smaller second-order risk for $d/c = 5$, 10 and 25, respectively. The degree of improvement decreases for larger values of $d/c$ as is anticipated since the expected number of stages of any reasonable test approaches 1 in this limit. Note that for each value of $d/c$, the expected number of stages of $\delta$ is larger than that of $\delta_{k^*}$, while the expected sample size is smaller. Thus, the way $\delta$ varies its stage sizes allows it to have more interim looks (stages), while keeping its overshoot, and hence, expected sample size, small. The expected number of stages and integrated risk of $\delta$ are close to their approximations $m_i^*$ and $\widetilde{r}_{c,d}$. Here the test $\delta$ was constructed from the samplers $\delta_{m,d/c}^o$. However, tests designed from the samplers $\delta_{m,z}^+$ also perform well in practice and behave almost identically to those constructed from the samplers $\delta_{m,d/c}^o$.

A natural question to ask is what values of $c, d$ should be used in practice if one is not comfortable specifying them as "costs"? The theory of the tests in Section 4 yields that $\log d^{-1}/I$ is an asymptotic approximation of the expected sample size and that $d$ is an asymptotic upper bound on the type I and II error probabilities. Hence, one could first choose $d$ to be the desired error probability or so that $\log d^{-1}/I$ is an acceptable expected sample size, and then choose $c$ so that $m_i^*$ is an acceptable expected number of stages, using (5.1).

## APPENDIX

**A.1. Proof of Theorems 2.1 and 2.2.** As mentioned above, the samplers $\delta_{m,h}^o$ and $\delta_{m,z}^+$ are extensions of Bartroff's [3] samplers for Brownian motion, and otherwise only differ slightly in their final stages. Moreover, Theorems 2.1 and 2.2 are extensions of Theorems 2.3 and 2.4 of Bartroff [3], requiring only two additional tools: first, justification for replacing the expected overshoot $E(\sum X_i - a; \sum X_i \geq a)$ by that of the normal distribution; second, bounds on the operating characteristics of the bold sampling $\hat{\delta}_n$ used in the final stages. With these two tools, the proofs of the corresponding theorems in Bartroff [3] can be followed almost exactly. We therefore state and prove these two needed tools here as Lemmas A.1 and A.2 and refer the reader to Bartroff [3] for the rest of the proof of Theorems 2.1 and 2.2. We also state without proof the auxiliary Lemma A.3 needed in the sequel, which is a simple extension of Lemma 2.4 of Bartroff [3] in the same manner.

Recall that $\psi^+(z) = \phi(z) - z\Phi(-z) = \int_z^\infty \Phi(-x)\,dx$. If $X_i$ are i.i.d. $N(\mu, 1)$ and $\Sigma_n = \sum_{i=1}^n X_i$, then

$$E(\Sigma_n - x; \Sigma_n \geq x) = \int_x^\infty P(\Sigma_n > y)\,dy = \sqrt{n} \cdot \psi^+\left(\frac{x - n\mu}{\sqrt{n}}\right).$$

Lemma A.1 shows that these two quantities are asymptotically equivalent in a certain range even when the $X_i$ are not normal, given that one of (2.3)–(2.5) holds.

LEMMA A.1.   *Let the $X_i$ be i.i.d. and satisfy one of (2.3)–(2.5). Let $a_n$ be a sequence such that*

$$\lim_{n \to \infty} \frac{a_n - n\mu}{\sqrt{n}} \in (-\infty, \infty) \quad or \quad \sqrt{(2-\varepsilon)\log n} \geq \frac{a_n - n\mu}{\sqrt{n}} \to \infty$$

*for some $\varepsilon \in (0,1)$ as $n \to \infty$. Then as $n \to \infty$,*

(A.3) $$P(\Sigma_n \geq a_n) \sim 1 - \Phi\left(\frac{a_n - n\mu}{\sqrt{n}}\right),$$

(A.4) $$E(\Sigma_n - a_n; \Sigma_n \geq a_n) \sim \sqrt{n} \cdot \psi^+\left(\frac{a_n - n\mu}{\sqrt{n}}\right).$$

PROOF.   Let $T_n = (\Sigma_n - n\mu)/\sqrt{n}$ and $b_n = (a_n - n\mu)/\sqrt{n}$. Assume that $b_n \to \infty$; otherwise, (A.3) holds by the central limit theorem. If (2.3) or (2.4) holds, then Theorem 4.6 of Hall [10] shows that $|P(T_n \geq x) - \Phi(-x)| = O(1/n)$ uniformly in $x$. Then

$$\left|\frac{P(T_n \geq b_n)}{\Phi(-b_n)} - 1\right| = \frac{O(1/n)}{\Phi(-b_n)} \leq \frac{O(1/n)}{\Phi(-\sqrt{(2-\varepsilon)\log n})}$$

$$= \frac{O(1/n)}{n^{-(2-\varepsilon)/2}/\sqrt{\log n}} = O(n^{-\varepsilon/2}\sqrt{\log n}) = o(1).$$

If (2.5) holds, then (A.3) holds by Cramér's theorem (e.g., see Feller [9], Theorem XVI.7.1).

Since $E(\Sigma_n - a_n; \Sigma_n > a_n) = \sqrt{n}\int_{b_n}^\infty P(T_n > x)\,dx$, to establish (A.4) it suffices to show that $\int_{b_n}^\infty P(T_n > x)\,dx \sim \psi^+(b_n)$. First assume that $b_n \to \infty$ such that $b_n \leq \sqrt{(2-\varepsilon)\log n}$. Choose $c_n \to \infty$ such that $b_n + \varepsilon' \leq c_n \leq \sqrt{(2-\varepsilon'')\log n}$, some $\varepsilon', \varepsilon'' > 0$. Then

(A.5) $$\frac{\phi(c_n)}{\psi^+(b_n)} \sim b_n^2 \frac{\phi(c_n)}{\phi(b_n)} \leq b_n^2 e^{-\varepsilon' c_n} \to 0$$

since $\psi^+(x) \sim \phi(x)/x^2$ as $x \to \infty$. Write $\int_{b_n}^\infty = \int_{b_n}^{c_n} + \int_{c_n}^\infty$. By (A.3),

(A.6) $$\int_{b_n}^{c_n} P(T_n > x)\,dx \sim \int_{b_n}^{c_n} \Phi(-x)\,dx = \psi^+(b_n) - \psi^+(c_n) \sim \psi^+(b_n)$$

since $\psi^+(c_n) \leq \phi(c_n) = o(\psi^+(b_n))$. For the other term,

$$\int_{c_n}^\infty P(T_n > x)\,dx = E(T_n; T_n > c_n) - c_n P(T_n > c_n)$$

by integration by parts and $c_n P(T_n > c_n) \sim c_n \Phi(-c_n) = o(\psi^+(b_n))$ by Mills' ratio and (A.5). By Schwarz's inequality, the other piece is

$$E(T_n; T_n > c_n) \leq \sqrt{ET_n^2 \cdot E1\{T_n > c_n\}^2}$$

$$= \sqrt{1 \cdot P(T_n > c_n)} \sim \sqrt{\Phi(-c_n)} = o(\psi^+(b_n))$$

by an argument like (A.5). These last two estimates give $\int_{c_n}^\infty P(T_n > x)\,dx = o(\psi^+(b_n))$, which with (A.6) gives $\int_{b_n}^\infty P(T_n > x)\,dx \sim \psi^+(b_n)$.

If $b_n \to b \in (-\infty, \infty)$, there is $T'_n$ with the same distribution as $T_n$ such that $T'_n \to Z \sim N(0,1)$ a.s. by weak convergence, and hence also in $L^1$ by uniform integrability (e.g., see Durrett [6], Theorems 2.1 and 5.2). Thus,

$$\int_{b_n}^\infty P(T_n > x)\,dx = E(T'_n - b_n; T'_n \geq b_n) \to E(Z - b; Z \geq b) = \psi^+(b). \quad \square$$

LEMMA A.2. *Let $n(x)$ be a positive integer-valued function and let $z(x) = (x - \mu n(x))/\sqrt{n(x)}$. If $n(x)$ is such that $z(x) \to -\infty$ and*

$$|z(x)| \leq [\sqrt{(2-\varepsilon)\log n(x)} \wedge \sqrt{x}] \tag{A.7}$$

*for some $\varepsilon \in (0,1)$ as $x \to \infty$, then $\hat{\delta}_{n(a)}(a) = (N, M)$ satisfies*

$$EN \leq a/\mu + O(|z(a)|\sqrt{a}) \tag{A.8}$$

*and $EM \to 1$ as $a \to \infty$.*

PROOF. Denote $n = n(a)$, $n_2 = \lceil n^{1/2} \rceil$ and $z = z(a)$. Suppose $y > 0$. It is well known from sequential theory that

$$E(M - 1|a - S_1 = y) \leq y/(\mu n_2) + O(1) = y/(\mu\sqrt{n}) + O(1)$$

as $n \to \infty$ uniformly in $y$. Thus,

$$\begin{aligned}EM - 1 &= E(M - 1; S_1 < a)\\ &\leq (\sqrt{n}\cdot\mu)^{-1} E(a - S_1; S_1 < a) + O(1) P(S_1 < a)\end{aligned} \tag{A.9}$$

and $(-a + \mu n)/\sqrt{n} = |z| \leq \sqrt{(2-\varepsilon)\log n}$, so by Lemma A.1,

$$E(a - S_1; S_1 < a) \sim \sqrt{n} \cdot \psi^+(|z|) \sim \sqrt{n} \cdot \frac{\phi(z)}{z^2}. \tag{A.10}$$

Also, since $P(S_1 < a) = P(\frac{S_1 - \mu n}{\sqrt{n}} < z) \to 0$, (A.9) becomes $EM \to 1$. To show that (A.8) holds, write $EN = n + n_2 \cdot E(M - 1) = n + o(\sqrt{n})$. We have $n = a/\mu + O(|z|\sqrt{a})$ by (A.7), so

$$EN = a/\mu + O(|z|\sqrt{a}) + o(\sqrt{a}) = a/\mu + O(|z|\sqrt{a}). \quad \square$$

LEMMA A.3. *If $h \in \mathcal{B}_m$ and $\delta$ is any sampler such that $R_h(\delta) = O(h(a))$, then for any $\varepsilon > 0$ and $0 \leq k < m$, as $a \to \infty$,*

$$P(a - S_k \geq (1 - \varepsilon)(1/\mu)^{1 - (1/2)^k} F_h^{(k)}(a)) \to 1.$$

**A.2. Proof of Theorem 3.1.** Let $a = \sigma^{-1}\log d^{-1}$ and $\{(c,d)\} \in \mathcal{B}_m^o$ so that $u_m(z^*) = m$. Let $a^* = \log d^{-1} + o(1)$ be that given by Lemma A.4 below. Then

$$h_k(\sigma^{-1}a^*) \sim \sigma^{-(1/2)^k} h_k(a^*) \propto h_k(\log d^{-1} + o(1)) \sim h_k(\log d^{-1})$$

since $(d/dx)h_k(x)$ is bounded for large $x$; thus,

(A.11) $$h_m(\sigma^{-1}a^*) \ll d/c \ll h_{m-1}(\sigma^{-1}a^*)$$

as a consequence of $\{(c,d)\} \in \mathcal{B}_m^o$. Let $\delta^* = (N^*, M^*)$ denote a Bayes power one test. By Lemma A.4, we know that $\sum_{i=1}^{N^*} Y_i = \sigma^{-1}\log l_{N^*} \geq \sigma^{-1}a^*$, so $\delta^*$ is a multistage sampler with boundary $\sigma^{-1}a^*$. Theorem 2.1 gives

$$\begin{aligned}
R_{c,d}^* &\geq cE_0 N^* + dE_0 M^* \\
&= c[E_0(N^* - a^*/I_0) + (d/c)E_0 M^*] + ca^*/I_0 \\
\text{(A.12)} \quad &\geq c[m(d/c) + o(d/c)] + cI_0^{-1}(\log d^{-1} + o(1)) \\
&= cI_0^{-1}\log d^{-1} + d \cdot m + o(d) \\
&= cI_0^{-1}\log d^{-1} + d \cdot u_m(z^*) + o(d).
\end{aligned}$$

Also by the $\mathcal{B}_m^o$ case of Theorem 2.1, for $\delta = (N, M)$,

(A.13) $$E_0(N - I_0^{-1}\log d^{-1}) + (d/c)E_0 M \leq m(d/c) + o(d/c).$$

Then

$$\begin{aligned}
R_{c,d}(\delta) &- P_1(N < \infty) \\
&= c[E_0(N - I_0^{-1}\log d^{-1}) + (d/c)E_0 M] + cI_0^{-1}\log d^{-1} \\
\text{(A.14)} \quad &\leq c[m(d/c) + o(d/c)] + cI_0^{-1}\log d^{-1} \quad \text{[by (A.13)]} \\
&= cI_0^{-1}\log d^{-1} + d \cdot m + o(d) \\
&= cI_0^{-1}\log d^{-1} + d \cdot u_m(z^*) + o(d),
\end{aligned}$$

so it suffices to show that $P_1(N < \infty) = o(d)$. The right-hand side of (A.13) is $O(d/c)$, so by Lemma A.3 (with $\sigma^{-1}I_0$ in place of $\mu$),

$$P_0(a - S_{m-1} \geq (1/2)(\sigma^{-1}I_0)^{-1+(1/2)^{m-1}} F_{d/c}^{(m-1)}(a)) \to 1$$

as $c, d \to 0$. On the above event $U$,

$$a - S_{m-1} \geq (1/2)(\sigma^{-1}I)^{-1+(1/2)^{m-1}} F_{d/c}^{(m-1)}(a) \geq \eta h_m(a)^2$$

for some $\eta > 0$ by Lemma 2.5 of Bartroff [3]. On $U$, the $m$th stage of $\delta = \delta_{m,d/c}^o(a)$ begins bold sampling. Letting $\rho_m = [(S_m - S_{m-1}) - \sigma^{-1}I_0 N_m]/\sqrt{N_m}$,

$$P_0(S_m \geq a + \sqrt{h_m(a)}|U) = P_0\left(\rho_m \geq \frac{a - S_{m-1} - \sigma^{-1}I_0 N_m}{\sqrt{N_m}} + \sqrt{\frac{h_m(a)}{N_m}}\bigg|U\right) \to 1$$

if $h_m(a) \ll N_m$ on $U$, since $(a - S_{m-1} - \sigma^{-1}I_0 N_m)/\sqrt{N_m} \to -\infty$ by definition of $\delta^o_{m,d/c}(a)$. This holds since

$$(\text{A.15}) \qquad N_m \geq \frac{a - S_{m-1}}{\sigma^{-1}I_0} \geq \frac{\eta h_m(a)^2}{\sigma^{-1}I_0} \gg h_m(a)$$

on $U$. Let $V = U \cap \{S_m \geq a + \sqrt{h_m(a)}\}$ so that $P_0(V) \to 1$. Using Wald's likelihood ratio identity, the relation $l_n = \exp(\sigma \sum_1^n Y_i)$, and letting primes denote complements,

$$\begin{aligned}
P_1(N < \infty) &= E_0(l_N^{-1}; N < \infty) \leq E_0 l_N^{-1} \\
&= E_0[\exp(-\sigma S_m); V] + E_0[\exp(-\sigma S_M); V'] \\
&\leq \exp(-\log d^{-1} - \sigma\sqrt{h_m(a)}) + E_0[\exp(-\log d^{-1}); V'] \\
&= d \cdot \exp(-\sigma\sqrt{h_m(a)}) + d \cdot P_0(V') \\
&= d \cdot o(1) + d \cdot o(1) = o(d),
\end{aligned}$$

proving that (3.6) holds in the $\{(c,d)\} \in \mathcal{B}^o_m$ case.

Now let $\{(c,d)\} \in \mathcal{B}^+_m$. By using the corresponding $\mathcal{B}^+_m$ cases of the results used in the arguments leading to (A.12) and (A.14),

$$R^*_{c,d} \geq cI_0^{-1}\log d^{-1} + d \cdot u_m(z^*) + o(d) \geq R_{c,d}(\delta) - P_1(N < \infty),$$

so it again suffices to show that $P_1(N < \infty) = o(d)$. Let $U$ be as above and $W_1 = \{S_m \geq a + \sqrt{h_m(a)}\}$, $W_2 = \{S_m \leq a - \sqrt{h_m(a)}\}$, $W_3 = \{S_{m+1} \geq a + (h_m(a))^{1/5}\}$ and $W = (U \cap W_1) \cup (U \cap W_2 \cap W_3)$. We will show that $P_0(W) \to 1$ as $d \to 0$, which will allow us to say that the likelihood ratio is large enough at the end of the $m$th stage (on $W_1$) or at the end of the $(m+1)$st stage (on $W_3$) that $P_1(N < \infty) = o(d)$:

$$\begin{aligned}
P_0(U \cap W_1) &= P_0(W_1|U)P_0(U) \sim P_0(W_1|U) \\
&= P_0\left(\rho_m \geq \frac{a - S_{m-1} - \sigma^{-1}I_0 N_m}{\sqrt{N_m}} + \sqrt{\frac{h_m(a)}{N_m}}\bigg|U\right)
\end{aligned}$$

and $(a - S_{m-1} - \sigma^{-1}I_0 N_m)/\sqrt{N_m} \to z^*$ on $U$ by definition of $\delta^+_{m,z^*}(a)$. Then $P_0(U \cap W_1) \to 1 - \Phi(z^*)$ by the central limit theorem if $\sqrt{h_m(a)} \ll \sqrt{N_m}$ on $U$, which holds by (A.15). Next, write

$$\begin{aligned}
P_0(U \cap W_2 \cap W_3) &= P_0(U)P_0(W_2|U)P_0(W_3|U \cap W_2) \\
&\sim P_0(W_2|U)P_0(W_3|U \cap W_2).
\end{aligned}$$

We have $P(W_2|U) \to \Phi(z^*)$ by an argument like that above. Also,

$$P_0(W_3|U \cap W_2) = P_0\left(\rho_{m+1} \geq \frac{a - S_m - \sigma^{-1}I_0 N_{m+1}}{\sqrt{N_{m+1}}} + \frac{(h_m(a))^{1/5}}{\sqrt{N_{m+1}}}\bigg|U \cap W_2\right),$$

which approaches 1 since

$$\sqrt{N_{m+1}} \geq \sqrt{\frac{a - S_m}{\sigma^{-1}I_0}} \geq \frac{(h_m(a))^{1/4}}{\sqrt{\sigma^{-1}I_0}} \gg (h_m(a))^{1/5}$$

and $(a - S_m - \sigma^{-1}I_0 N_{m+1})/\sqrt{N_{m+1}} \to -\infty$ on $U \cap W_2$ by definition of $\delta^+_{m,z^*}(a)$. Combining these, we have $P_0(U \cap W_2 \cap W_3) \to \Phi(z^*)$ and, hence,

$$P_0(W) = P_0(U \cap W_1) + P_0(U \cap W_2 \cap W_3) \to 1 - \Phi(z^*) + \Phi(z^*) = 1.$$

Note that on $W$, $S_M - a \geq \sqrt{h_m(a)} \wedge (h_m(a))^{1/5} = (h_m(a))^{1/5}$, so

$$\begin{aligned}
P_1(N < \infty) &= E_0(l_N^{-1}; N < \infty) \leq E_0 l_N^{-1} \\
&= E_0[\exp(-\sigma S_M); W] + E_0[\exp(-\sigma S_M); W'] \\
&\leq \exp(-\log d^{-1} - \sigma(h_m(a))^{1/5}) + E_0[\exp(-\log d^{-1}); W'] \\
&= d \cdot \exp(-\sigma(h_m(a))^{1/5}) + d \cdot P_0(W') \\
&= d \cdot o(1) + d \cdot o(1) = o(d),
\end{aligned}$$

finishing the proof.

LEMMA A.4. *There exists $a^* = \log d^{-1} + o(1)$ such that $\log l_{N^*} \geq a^*$.*

PROOF. Suppose that a Bayes procedure has sampled $X_1, \ldots, X_n$ in $m$ stages. By the Bayes property, $\delta^*$ will stop at this point only if the stopping risk is no greater than the continuation risk, that is, only if

(A.16) $$l_n^{-1} \leq \rho(c, d, l_n^{-1}),$$

where

$$\rho(u, v, w) = \inf_{(N,M): N \geq 1} \{E_0(uN + vM) + wP_1(N < \infty)\}.$$

Multiplication of (A.16) by $l_n$ yields $1 \leq \rho(l_n c, l_n d, 1)$; hence, we consider the function $\rho(t) = \rho(tc, td, 1)$ for $t > 0$, and note that (A.16) implies that $\rho(l_{N^*}) \geq 1$. The function $\rho(t)$ is the infimum of a set of lines, each of slope at least $c + d$ by virtue of the restriction on the infimum. Thus, $\rho(t)$ is continuous, strictly increasing and satisfies $\rho(t) \geq t(c + d)$, so that

(A.17) $$\rho(t) \geq 1 \quad \text{when } t \geq (c + d)^{-1}.$$

If $(N', M')$ is the procedure that samples with constant stage size one (i.e., fully-sequential sampling) and an appropriately chosen boundary, then it is well known (e.g., see Lorden [16]) that $P_1(N' < \infty) < 1$ and $E_0 N' = E_0 M' < \infty$, and hence, $\rho(t) \leq t(c+d)E_0 N' + P_1(N' < \infty) < 1$ for sufficiently small $t$. This and (A.17) imply that there is a unique number $e^{a^*}$ such that $\rho(e^{a^*}) = 1$. Then $\log l_{N^*} = \log \rho^{-1}(\rho(l_{N^*})) \geq \log \rho^{-1}(1) = a^*$. To show that $a^* = \log d^{-1} + o(1)$, let $Y_i$ be as in (3.4) and $\delta_{1,h}^o(a) = (N, M)$, the multistage sampler described in Section 2 with $h(a) = a^{3/4}$ and $a = \sigma^{-1}\log(d/c)$. Since $\sqrt{a} \ll h(a) \ll a$, by Lemma A.2, $E_0 N - a(\sigma^{-1}I_0)^{-1} = o(h(a))$ and $E_0 M \to 1$. Also, $l_N^{-1} = \exp[-\sigma(Y_1 + \cdots + Y_N)] \leq \exp[-\sigma a] = c/d$, so that

$$\rho(t) \leq E_0[t(cN + dM) + l_N^{-1}1\{N < \infty\}]$$
$$\leq tc[E_0(N - a(\sigma^{-1}I_0)^{-1}) + a(\sigma^{-1}I_0)^{-1} + (d/c)E_0 M] + E_0 l_N^{-1}$$
$$\leq tc[o(h(a)) + a(\sigma^{-1}I_0)^{-1} + (d/c)(1 + o(1))] + c/d$$
$$= tc[o(d/c) + d/c(1 + o(1))] + c/d = td(1 + o(1)) + c/d.$$

This implies that $\rho(t) \leq 1$ when $t \leq d^{-1}(1 + o(1))$; hence,

$$a^* = \log \rho^{-1}(1) \geq \log \rho^{-1}(\rho(d^{-1}(1 + o(1)))) = \log d^{-1} + o(1).$$

On the other hand,

$$a^* = \log \rho^{-1}(1) \leq \log \rho^{-1}(\rho([c + d]^{-1})) \qquad \text{[by (A.17)]}$$
$$= \log(c + d)^{-1} = \log d^{-1} + o(1)$$

since $d/c \to \infty$, establishing $a^* = \log d^{-1} + o(1)$. $\square$

### A.3. Proof of Theorems 4.1 and 4.2.

LEMMA A.5. *Assume $\{(c,d)\} \in \mathcal{B}_m$ and let $z^* \in [-\infty, \infty)$ be the unique solution of*

$$\frac{\phi(z^*)}{1 - \Phi(z^*)} = \lim_{c,d \to 0} \frac{\kappa_m(\sigma_0^{-1}I_0)h_m(\sigma_0^{-1}\log d^{-1})}{d/c}.$$

*Then*

(A.18) $\quad cE_0 N^* + dE_0 M^* + P_1(D^* = 0) \geq cI_0^{-1}\log d^{-1} + u_m(z^*)d - o(d)$

*as $c, d \to 0$.*

PROOF. We extend $\delta^*$ to a power one test of $f_0$ versus $f_1$ on the event $\{D^* = 1\}$. Let $N = M = \inf\{n \geq 1 : l_n \geq d^{-2}\}$ be fully-sequential sampling with likelihood ratio boundary $d^{-2}$. Define $N' = N^* + N \cdot 1\{D^* = 1\}$ and $M' = M^* + M \cdot 1\{D^* = 1\}$, the power one test that coincides with $\delta^*$ on $\{D^* = 0\}$ but continues with the power one test $(N, M)$ on $\{D^* = 1\}$. Since $\{N' < \infty\} = \{D^* = 0\} \cup \{D^* = 1, N < \infty\}$, we have

$$cE_0 N^* + dE_0 M^* + P_1(D^* = 0)$$
$$= c[E_0 N' - E_0(N; D^* = 1)] + d[E_0 M' - E_0(M; D^* = 1)]$$

$$+ P_1(N' < \infty) - P_1(D^* = 1, N < \infty)$$
$$= [cE_0 N' + dE_0 M' + P_1(N' < \infty)]$$
$$\quad - [cE_0(N; D^* = 1) + dE_0(M; D^* = 1) + P_1(D^* = 1, N < \infty)]$$
$$= R_1 - R_2.$$

By Theorem 3.1, $R_1 \geq cI^{-1} \log d^{-1} + u_m(z^*)d + o(d)$, so to show that (A.18) holds, it suffices to show that $R_2 = o(d)$. Write

$$\text{(A.19)} \quad \begin{aligned} R_2 &\leq [cE_0(N|D^* = 1) + dE_0(M|D^* = 1)]P_0(D^* = 1) \\ &\quad + P_1(N < \infty | D^* = 1). \end{aligned}$$

It is well known that

$$E_0(N|D^* = 1) = E_0(M|D^* = 1) = I_0^{-1} \log d^{-2} + O(1) = O(\log d^{-1}).$$

We will show below that there is a $K < \infty$ such that

$$\text{(A.20)} \quad l_{N^*} \leq Kd \quad \text{on } \{D^* = 1\}.$$

Using this and Wald's likelihood identity,

$$P_0(D^* = 1) = E_1(l_{N^*}; D^* = 1, N^* < \infty) \leq Kd = O(d).$$

Combining these two estimates gives

$$\text{(A.21)} \quad \begin{aligned} [cE_0(N|D^* &= 1) + dE_0(M|D^* = 1)]P_0(D^* = 1) \\ &= [c \cdot O(\log d^{-1}) + d \cdot O(\log d^{-1})]O(d) \\ &= O(d^2 \log d^{-1}). \end{aligned}$$

By definition of $(N, M)$,

$$\begin{aligned} P_1(N < \infty | D^* = 1) &= E_0(l_N^{-1} 1\{N < \infty\} | N > 0) \\ &\leq E_0(d^2 1\{N < \infty\} | N > 0) \leq d^2. \end{aligned}$$

Plugging this and (A.21) into (A.19) gives $R_2 \leq O(d^2 \log d^{-1}) + d^2 = o(d)$.

To verify (A.20), write the posterior risk $r_{ik}$ of rejecting $H_i$ after the $k$th stage as

$$\text{(A.22)} \quad r_{0k} = \frac{w_0 \pi_0 l_{N^* k}}{\pi_0 l_{N^* k} + \pi_1}, \qquad r_{1k} = \frac{w_1 \pi_1}{\pi_0 l_{N^* k} + \pi_1},$$

and let $r_k = r_{0k} \wedge r_{1k}$, the stopping risk after the $k$th stage. A Bayes test stops sampling if the stopping risk is less than all possible continuation risks. One possible continuation is fully-sequential sampling. By Lemma 2 of Lorden [14] there is a constant $K^* < \infty$ such that a Bayes procedure can only stop when the continuation risk of fully-sequential sampling is less than $K^*$ times the cost per observation, $c + d$ in this case. Thus, $r_{M^*} \leq K^*(c+d) \leq 2K^*d$, meaning $r_{0M^*} \leq 2K^*d$ or $r_{1M^*} \leq 2K^*d$. If $r_{0M^*} \leq 2K^*d$, then by the first relation in (A.22) and some simple algebra,

$$l_{N^*} \leq \frac{\pi_1 \cdot 2K^*d}{\pi_0(w_0 - 2K^*d)} \leq \frac{4\pi_1 K^*}{\pi_0 w_0} d$$

for small enough $d$. Clearly, $r_{0M^*} < r_{1M^*}$ in this case, so we can be sure $D^* = 1$. Otherwise, $r_{1M^*} \leq 2K^*d < r_{0M^*}$ for small $d$, so $D^* = 0$. □

PROOF OF THEOREM 4.1. Let $I = I_0 = I_1$ and $\sigma = \sigma_0 = \sigma_1$. Rearranging terms,

(A.23) $$r^*_{c,d} = \sum_{i=0}^{1} \pi_{1-i} w_{1-i} [c_i E_i N^* + d_i E_i M^* + P_{1-i}(D^* = i)],$$

where $c_i = c\pi_i/(\pi_{1-i} w_{1-i})$ and $d_i = d\pi_i/(\pi_{1-i} w_{1-i})$. It is simple to verify that $\{(c_i, d_i)\} \in \mathcal{B}_m$ and

$$\lim_{c,d \to 0} \frac{\kappa_m(\sigma^{-1} I) h_m(\sigma^{-1} \log d_i^{-1})}{d_i/c_i} = \frac{\phi(z^*)}{1 - \Phi(z^*)}.$$

By Lemma A.5,

$$c_i E_i N^* + d_i E_i M^* + P_{1-i}(D^* = i) \geq c_i I^{-1} \log d_i^{-1} + u_m(z^*) d_i + o(d_i),$$

and plugging this into (A.23) gives

$$r^*_{c,d} \geq \sum_{i=0}^{1} \pi_{1-i} w_{1-i} [c_i I^{-1} \log d^{-1} + u_m(z^*) d_i + o(d_i)]$$

$$= \sum_{i=0}^{1} \pi_i [c I^{-1} \log d^{-1} + u_m(z^*) d + o(d)]$$

$$= c I^{-1} \log d^{-1} + u_m(z^*) d + o(d),$$

establishing (4.1). For an event $A$, denote

(A.24) $$r_{c,d}(\delta; A) = \sum_{i=0}^{1} \pi_i [c E_i(N; A) + d E_i(M; A) + w_i P_i(D = 1 - i, A)].$$

Obviously $r_{c,d}(\delta; A) + r_{c,d}(\delta; A') = r_{c,d}(\delta)$. Let

$$A_0 = \{|\log l^1 - IN_1| \leq \sigma \sqrt{\beta N_1 \log N_1}\},$$
$$A_1 = \{|-\log l^1 - (-IN_1)| \leq \sigma \sqrt{\beta N_1 \log N_1}\},$$

where $\beta = 2 - (1/2)^{m-1}$. The following bounds are proved below:

(A.25) $$cE_0(N; A_0) \leq cE_0 N^{(0)} + o(d),$$

(A.26) $$dE_0(M; A_0) \leq dE_0 M^{(0)} + o(d),$$

(A.27) $$P_0(D = 1, A_0) = o(d),$$

(A.28) $$cE_1(N; A_0) = o(d),$$

(A.29) $$dE_1(M; A_0) = o(d),$$

(A.30) $$P_1(D = 0, A_0) \leq P_1(N^{(0)} < \infty) + o(d).$$

Using these bounds,

$$r_{c,d}(\delta; A_0)$$
$$= \sum_{i=0}^{1} \pi_i [cE_i(N; A_0) + dE_i(M; A_0) + w_i P_i(D = 1 - i, A_0)]$$
$$\leq \pi_0 [cE_0 N^{(0)} + dE_0 M^{(0)} + o(d)] + \pi_1 [w_1 P_1(N^{(0)} < \infty) + o(d)]$$
$$= \pi_1 w_1 [c_0 E_0 N^{(0)} + d_0 E_0 M^{(0)} + P_1(N^{(0)} < \infty)] + o(d)$$
$$= \pi_1 w_1 [c_0 I^{-1} \log d_0^{-1} + u_m(z^*) d_0 + o(d_0)] + o(d) \qquad \text{(by Theorem 3.1)}$$
$$= \pi_0 [cI^{-1} \log d^{-1} + u_m(z^*) d] + o(d)$$

and the same argument with the indices reversed yields

(A.31) $$r_{c,d}(\delta; A_1) \leq \pi_1 [cI^{-1} \log d^{-1} + u_m(z^*) d] + o(d).$$

Now we consider $r_{c,d}(\delta; A_0' \cap A_1')$. Let $A = A_0' \cap A_1'$. The bounds

(A.32) $$cE_0(N; A) = o(d),$$

(A.33) $$dE_0(M; A) = o(d),$$

(A.34) $$P_0(D = 1, A) = o(d),$$

are also proved below. These bounds give $r_{c,d}(\delta; A) = o(d)$. Combining this with (A.31) gives

$$r_{c,d}(\delta) = r_{c,d}(\delta; A_0) + r_{c,d}(\delta; A_1) + r_{c,d}(\delta; A)$$
$$\leq \sum_{i=0}^{1} \pi_i [cI^{-1} \log d^{-1} + u_m(z^*) d] + o(d)$$
$$= cI^{-1} \log d^{-1} + d \cdot u_m(z^*) + o(d),$$

establishing (4.2). All that remains is to verify the bounds (A.25)–(A.30) and (A.32)–(A.34).

Let $Y_j = Y_j^{(0)}$ and $n(x,z)$ be as in Section 2 with $\sigma^{-1}I$ in place of $\mu$. We begin by proving the crude bound

(A.35) $\quad E_i(N|U) = O(\log d^{-1}) \quad$ for any $U$ such that $E_i(M|U) = O(1)$,

$i = 0, 1$. If $\{(c,d)\} \in \mathcal{B}_m^o$, the $m$th stage begins bold sampling, in which each stage is bounded by $\max_{x \in \{x_+, x_-\}} n(x, -[\sqrt{d/c}/x^{1/4} \wedge \sqrt{(3/2)\log(x+1)}])$, where $x_+ = \sigma^{-1}\log d^{-1} - \sum Y_j$ and $x_- = \sum Y_j - (-\sigma^{-1}\log d^{-1})$. If $\delta$ does not stop at the end of a stage, then it must be that $|\sum Y_j| < \sigma^{-1}\log d^{-1}$; hence, $x_+$ and $x_-$ are both bounded above by $2\sigma^{-1}\log d^{-1}$. Then the sizes of stages $m, m+1, \ldots$ are all bounded by

$$n(x, -[\sqrt{d/c}/x^{1/4} \wedge (1/2)\sqrt{\log x}])|_{x=2\sigma^{-1}\log d^{-1}} = O(\log d^{-1}).$$

The sizes of the first $m-1$ stages are likewise bounded by

$$\max_{x \in \{x_+, x_-\}} n(x, \sqrt{(1-2^k)\log(x+1)}) \leq n(x, \sqrt{(1/2)\log(x+1)})|_{x=2\sigma^{-1}\log d^{-1}}$$
$$= O(\log d^{-1})$$

for some $k \geq 1$. Thus, the size of each stage of $\delta$ is uniformly $O(\log d^{-1})$ and therefore, $E_i(N|U) \leq O(\log d^{-1})E_i(M|U) = O(\log d^{-1})$. This holds if $\{(c,d)\} \in \mathcal{B}_m^+$ as well since then the $(m+1)$st stage begins bold sampling. Next let $B = \{\log l^k > -\log d^{-1}$ for all $1 \leq k \leq M\}$ and note that $\delta$ and $(N^{(0)}, M^{(0)})$ coincide on $A_0 \cap B$ since $\log l^1 \geq IN_1 - \sigma^{-1}\sqrt{\beta N_1 \log N_1} > 0$ for small $d$ on $A_0$ and $\log l^k$ never crosses the lower boundary $-\log d^{-1}$ on $B$. Clearly, $E_0(M|A_0 \cap B') = O(1)$, so using this crude bound and Wald's likelihood identity,

$$P_0(A_0 \cap B') \leq P_0(B') = E_1(l^M; B') \leq E_1(d; B') \leq d$$

and $E_0(N; A_0 \cap B) \leq E_0 N^{(0)}$ since $\delta$ and $(N^{(0)}, M^{(0)})$ coincide on $A_0 \cap B$, so that

$$cE_0(N; A_0) = cE_0(N; A_0 \cap B) + cE_0(N; A_0 \cap B')$$
$$\leq cE_0 N^{(0)} + c \cdot O(d \log d^{-1})$$
$$= cE_0 N^{(0)} + o(c) = cE_0 N^{(0)} + o(d),$$

which proves (A.25). Similarly, $E_0(M; A_0 \cap B) \leq E_0 M^{(0)}$ and $E_0(M|A_0 \cap B') = O(1)$, so that

$$dE_0(M; A_0) \leq dE_0(M; A_0 \cap B) + dE_0(M|A_0 \cap B')P_0(A_0 \cap B')$$
$$\leq dE_0 M^{(0)} + d \cdot O(1) \cdot d = dE_0 M^{(0)} + o(d),$$

proving (A.26). Letting $\gamma(d) = IN_1 - \sigma\sqrt{\beta N_1 \log N_1}$,

$$P_0(D = 1, A_0) \leq P_0(D = 1|A_0) = P_0(l^M \leq -\log d^{-1} | \log l^1 \geq \gamma(d))$$
$$\leq \exp[-(\log d^{-1} + \gamma(d))] = o(d),$$

proving (A.27), and a similar argument proves (A.34). Since $\gamma \sim IN_1 \sim \log d^{-1}$, we have

$$P_1(A_0) = E_0(l_1^{-1}; \log l_1 \geq \gamma(d)) \leq E_0(e^{-\gamma(d)}; \log l_1 \geq \gamma(d))$$
$$\leq e^{-\gamma(d)} \leq \exp[-(1/2)\log d^{-1}] = \sqrt{d}.$$

Also, $E_1(N|A_0) = O(\log d^{-1})$ by (A.35) so

$$cE_1(N; A_0) = cE_1(N|A_0)P_1(A_0) \leq c\sqrt{d} \cdot O(\log d^{-1}) = c \cdot o(1) = o(d),$$

proving (A.28). (A.29) holds since $E_1(M|A_0) = O(1)$ and $P_1(A_0) \to 0$ and similarly for (A.33). Since $\delta$ and $(N^{(0)}, M^{(0)})$ coincide on $A_0 \cap B$,

$$P_1(D = 0, A_0 \cap B) = P_1(N^{(0)} < \infty, A_0 \cap B) \leq P_1(N^{(0)} < \infty).$$

Also,

$$P_1(D = 0, A_0 \cap B') = E_0[(l^M)^{-1}; D = 0, A_0 \cap B']$$
$$\leq E_0[d; D = 0, A_0 \cap B']$$
$$\leq dP_0(B') = o(d)$$

since clearly $P_0(B') \to 0$. Combining these two gives

$$P_1(D = 0; A_0) = P_1(D = 0; A_0 \cap B) + P_1(D = 0; A_0 \cap B')$$
$$\leq P_1(N^{(0)} < \infty) + o(d),$$

proving (A.30). Now

$$P_0(A) \leq P_0(A_0') = P_0(\log l^1 < \gamma(d)) = P_0\left(\frac{-\log l^1 + IN_1}{\sigma\sqrt{N_1}} > \frac{IN_1 - \gamma(d)}{\sigma\sqrt{N_1}}\right)$$

and $(IN_1 - \gamma(d))/(\sigma\sqrt{N_1}) = \sqrt{\beta \log N_1}$, so by (A.3),

$$P_0(A) \lesssim \Phi(-\sqrt{\beta \log N_1}) \sim \frac{N_1^{-\beta/2}}{\sqrt{\log N_1}} = o((\log d^{-1})^{-\beta/2})$$

since $IN_1 \sim \log d^{-1}$. Then since $d/c = O(h_m(\log d^{-1}))$,

$$cE_0(N; A) = cE_0(N|A)P_0(A) = c \cdot O(\log d^{-1}) \cdot o((\log d^{-1})^{-\beta/2})$$

(A.36)
$$= o(d) \cdot \frac{(\log d^{-1})^{1-\beta/2}}{d/c} = o(d) \cdot \frac{(\log d^{-1})^{(1/2)^m}}{h_m(\log d^{-1})}$$

$$= o(d) \cdot \frac{(\log d^{-1})^{(1/2)^m}}{(\log d^{-1})^{(1/2)^m}} = o(d),$$

proving (A.32) and finishing the proof. □

PROOF OF THEOREM 4.2. Let $T = \{t > 0 : |\log t - I_0 N_1| \leq \sigma_0 \sqrt{\beta N_1 \log N_1}\}$, where $\beta = 2 - (1/2)^{m-1}$ and $A_0 = \{l^1 \in T\}$. Let $\delta = (N, M, D)$ and $\delta_0(l^1 c, l^1 d, z_0^*) = (N^{(0)}, M^{(0)})$. We will use the $r_{c,d}(\delta; A)$ notation as in (A.24). Clearly, $l^1 \geq 1$ on $A_0$ for sufficiently small $c, d$, so that $\delta$ will switch and continue sampling according to $\delta_0(l^1 c, l^1 d, z_0^*)$ after the first stage; hence,

$$\text{(A.37)} \qquad N \leq N^{(0)} + N_1 \quad \text{and} \quad M \leq M^{(0)} + 1.$$

Also note that

$$\text{(A.38)} \qquad \{D = 0\} \cap A_0 \subseteq \{N^{(0)} < \infty\},$$

since on $\{D=0\} \cap A_0$ the likelihood ratio will cross the boundary $d^{-1}$, which is equivalent to the stopping rule of $\delta_0(l^1 c, l^1 d, z_0^*)$, as discussed above. Using the bounds (A.27)–(A.29),

$$\text{(A.39)} \quad \begin{aligned} r_{c,d}(\delta; A_0) &= \pi_0 c E_0(N; A_0) + \pi_0 d E(M; A_0) \\ &\quad + \pi_1 w_1 P_1(D = 0, A_0) + o(d) \\ &= E_0[\pi_0 c N + \pi_0 d M + \pi_1 w_1 (l^M)^{-1} \cdot 1\{D=0\}; A_0] + o(d) \\ &\leq E_0[\pi_0 c N^{(0)} + \pi_0 d M^{(0)} \\ &\quad + \pi_1 w_1 (l^1)^{-1} (l^M/l^1)^{-1} \cdot 1\{N^{(0)} < \infty\}; A_0] \\ &\quad + \pi_0 c N_1 + \pi_0 d + o(d) \qquad [\text{by (A.37) and (A.38)}] \\ &= E_0[\varphi(l^1); l^1 \in T] + \pi_0 c N_1 + \pi_0 d + o(d), \end{aligned}$$

where

$$\varphi(t) = c E_0[\pi_0(c N^{(0)} + d M^{(0)}) | l^1 = t] + \pi_1 w_1 t^{-1} P_1(N^{(0)} < \infty | l^1 = t).$$

By rearranging terms,

$$\text{(A.40)} \quad \begin{aligned} \varphi(t) &= \pi_0 t^{-1} \{E_0[(tc) N^{(0)} + (td) M^{(0)} | l^1 = t] + P_1(N^{(0)} < \infty | l^1 = t)\} \\ &\quad + (\pi_1 w_1 - \pi_0) t^{-1} P_1(N^{(0)} < \infty | l^1 = t) \\ &= \pi_0 t^{-1} R_{tc, td}(\delta(tc, td, z_0^*)) + (\pi_1 w_1 - \pi_0) t^{-1} P_1(N^{(0)} < \infty | l^1 = t). \end{aligned}$$

For any $t \in T$, $\log(td)^{-1} \sim (1 - I_0/I_1) \log d^{-1}$, which implies that

$$h_m(\log(td)^{-1}) \sim h_m((1 - I_0/I_1) \log d^{-1}) \sim (1 - I_0/I_1)^{(1/2)^m} h_m(\log d^{-1}),$$

and hence, $\{(tc, td)\} \in \mathcal{B}_m$ uniformly for $t \in T$. Moreover, by this last,

$$\text{(A.41)} \qquad \lim_{c,d \to 0} \frac{\kappa_m(\sigma_0^{-1} I_0) h_m(\sigma_0^{-1} \log(td)^{-1})}{(td)/(tc)} = \frac{\phi(z_0^*)}{1 - \Phi(z_0^*)}.$$

By Theorem 3.1,
$$R_{tc,td}(tc,td,z_0^*) \le tcI_0^{-1}\log(td)^{-1} + u_m(z_0^*)td + o(td)$$

and the proof of Theorem 3.1 shows that $P_1(N^{(0)} < \infty | l^1 = t) = o(td)$ uniformly for $t \in T$. Plugging these last two into (A.40),

$$\varphi(t) \le \pi_0 t^{-1}[tcI_0^{-1}\log(td)^{-1} + u_m(z_0^*)td + o(td)] + (\pi_1 w_1 - \pi_0)t^{-1}o(td)$$
$$= \pi_0[cI_0^{-1}\log d^{-1} + u_m(z_0^*)d] - \pi_0 cI_0^{-1}\log t + o(d),$$

uniformly on $T$, and, in turn, plugging this into (A.39) gives

(A.42)
$$r_{c,d}(\delta; A_0) \le \pi_0[cI_0^{-1}\log d^{-1} + (1 + u_m(z_0^*))d]$$
$$+ \pi_0 cI_0^{-1}[I_0 N_1 - E(\log l^1; l^1 \in T)] + o(d).$$

By repeating the argument leading to (A.36), we have $E_0(\log l^1; A_0) = o(d/c)$; hence, (A.42) becomes

(A.43) $\quad r_{c,d}(\delta; A_0) \le \pi_0[cI_0^{-1}\log d^{-1} + (1 + u_m(z_0^*))d] + o(d).$

Letting $A_1 = \{|\log(1/l^1) - I_1 N_1| \le \sigma_1\sqrt{\beta N_1 \log N_1}\}$ and repeating arguments in the proof of Theorem 4.1 give

$$r_{c,d}(\delta; A_1) \le \pi_1[cI_1^{-1}\log d^{-1} + d \cdot u_m(Q_1, \sigma_1^{-1} I_1)] + o(d)$$

and $r_{c,d}(\delta; A_0' \cap A_1') = o(d)$. Combining with (A.43) gives (4.8) with a "$\le$."

Next we show that (4.7) holds with a "$\ge$." Let $l^{*k} = l_{N^{*k}}$, $T^* = \{t > 0 : |\log t - I_0 N_1^*| \le \sigma_0\sqrt{\beta N_1^* \log N_1^*}\}$ and $A_0^* = \{l^{*1} \in T^*\}$. Let

$$r_i^* = \pi_i(cE_i N^* + dE_i M^*) + \pi_{1-i} w_{1-i} P_{1-i}(D^* = i), \qquad i = 0, 1.$$

Since $\delta^*$ follows its first stage with the optimal continuation $(\dot N^*, \dot M^*, \dot D^*)$, we can write

(A.44)
$$r_0^* = E_0[\pi_0(cN^* + dM^*) + \pi_1 w_1(l^{*M^*})^{-1}1\{D^* = 0\}]$$
$$= E_0[\pi_0(cE_0 \dot N^* + dE_0 \dot M^*) + \pi_1 w_1(l^{*1})^{-1}P_1(\dot D^* = 0)]$$
$$+ \pi_0(cN_1^* + d).$$

Define $\varphi^*(t) = \pi_1 w_1 t^{-1}\{E_0[c(t)\dot N^* + d(t)\dot M^* | l^{*1} = t] + P_1(\dot D^* = 0 | l^{*1} = t)\}$, where $c(t) = ct\pi_0/(\pi_1 w_1)$ and $d(t) = dt\pi_0/(\pi_1 w_1)$. It will be shown below that $N_1^* \sim I_1^{-1}\log d^{-1}$. Assuming this holds, the arguments leading to (A.41) show that it holds with $(tc, td)$ replaced by $(c(t), d(t))$. Then by Lemma A.5,

(A.45)
$$\varphi^*(t) \ge \pi_1 w_1 t^{-1}[c(t)I_0^{-1}\log d(t)^{-1} + u_m(z_0^*)d(t) + o(d(t))]$$
$$= \pi_0[cI_0^{-1}\log d^{-1} + u_m(z_0^*)d - \pi_0 cI_0^{-1}\log t + o(d)$$

uniformly for $t \in T^*$, and hence,

$$\begin{aligned}
r_0^* &= E_0[\varphi^*(l^{*1})] + \pi_0(cN_1^* + d) \\
&\geq E_0[\varphi^*(l^{*1}); A_0^*] + \pi_0(cN_1^* + d) \qquad \text{(since } \varphi^* \geq 0\text{)} \\
&\geq \pi_0[cI_0^{-1}\log d^{-1} + u_m(z_0^*)d]P_0(A_0^*) - \pi_0 cI_0^{-1}E_0[\log l^{*1}; A_0^*] \\
&\quad + \pi_0(cN_1^* + d) + o(d) \qquad \text{[by (A.45)]} \\
&\geq \pi_0[cI_0^{-1}\log d^{-1} + (1 + u_m(z_0^*))d] + o(d),
\end{aligned}$$
(A.46)

this last by the arguments leading to (A.36). A straightforward application of Lemma A.5 gives $r_1^* \geq \pi_1[cI_1^{-1}\log d^{-1} + u_m(z_1^*)d] + o(d)$, and adding these last two gives (4.7).

All that remains is to verify that $N_1^* \sim I_1^{-1}\log d^{-1}$. Suppose instead that

$$\underline{L} = \liminf_{c,d \to 0} \frac{N_1^*}{\log d^{-1}} < I_1^{-1}.$$
(A.47)

Then there is a sequence $\{(c,d)\}$ approaching $(0,0)$ on which the lim inf is achieved, and by repeating the above arguments on this sequence,

$$r_0^* \geq \pi_0[cI_0^{-1}\log d^{-1} + (1 + u_m(z_0'))d] + o(d),$$
(A.48)

where $z_0'$ is the unique solution of

$$\frac{\phi(z_0')}{1 - \Phi(z_0')} = \lim_{c,d \to 0} \frac{\kappa_m(\sigma_0^{-1}I_0)h_m(\sigma_0^{-1}(1 - I_0\underline{L})\log d^{-1})}{d/c}.$$

By writing

$$h_m(\sigma_0^{-1}(1 - I_0\underline{L})\log d^{-1}) = \left(\frac{1 - I_0\underline{L}}{1 - I_0/I_1}\right)^{(1/2)^m} \times h_m(\sigma_0^{-1}(1 - I_0/I_1)\log d^{-1}),$$

we have

$$\begin{aligned}
\frac{\phi(z_0')}{1 - \Phi(z_0')} &= \left(\frac{1 - I_0\underline{L}}{1 - I_0/I_1}\right)^{(1/2)^m} \\
&\quad \times \lim_{c,d \to 0} \frac{\kappa_m(\sigma_0^{-1}I_0)h_m(\sigma_0^{-1}(1 - I_0/I_1)\log d^{-1})}{d/c} \\
&= \left(\frac{1 - I_0\underline{L}}{1 - I_0/I_1}\right)^{(1/2)^m} \times \frac{\phi(z_0^*)}{1 - \Phi(z_0^*)} \geq \frac{\phi(z_0^*)}{1 - \Phi(z_0^*)}.
\end{aligned}$$

Hence, $z_0' \geq z_0^*$ since $z \mapsto \phi(z)/[1 - \Phi(z)]$ is increasing, so (A.48) becomes

$$r_0^* \geq \pi_0[cI_0^{-1}\log d^{-1} + (1 + u_m(z_0^*))d] + o(d),$$
(A.49)

since $u_m$ is strictly increasing. By reversing indices and repeating this argument, conditioning on $\{|\log(1/l^{*1}) - I_1 N_1^*| \leq \sigma_1\sqrt{\beta N_1^* \log N_1^*}\}$ instead of $A_0^*$, we obtain

$$r_1^* \geq \pi_1[cI_1\log d^{-1} + (1 + u_m(z_1'))d] + o(d).$$
(A.50)

Using (A.49), (A.50) and (4.7), we would then have

$$\begin{aligned}
r_{c,d}^* - r_{c,d}(\delta) &= r_0^* + r_1^* - r_{c,d}(\delta) \\
&\geq \pi_1 d[1 + u_m(z_1') - u_m(z_1^*)] + o(d) \\
&\geq \varepsilon d + o(d) > 0
\end{aligned}$$

for some $\varepsilon > 0$ and sufficiently small $c, d$ since $m \leq u_m < m + 1$. This obviously contradicts $r_{c,d}^* \leq r_{c,d}(\delta)$ so (A.47) cannot hold. On the other hand, if

$$\text{(A.51)} \qquad \eta = \limsup_{c,d \to 0} \frac{N_1^*}{\log d^{-1}} - I_1^{-1} > 0,$$

then again on a sequence $\{(c,d)\}$ approaching $(0,0)$, we would have

$$\begin{aligned}
r_{c,d}^* - r_{c,d}(\delta) &\geq \pi_0 c I_0^{-1} \log d^{-1} + \pi_1 c N_1^* - r_{c,d}(\delta) \qquad \text{(by Lemma A.5)} \\
&\geq \pi_0 c I_0^{-1} \log d^{-1} + \pi_1 c(\eta + I_1^{-1}) \log d^{-1}(1 + o(1)) \\
&\quad - [(\pi_0/I_0 + \pi_1/I_1) c \log d^{-1} + O(d)] \qquad \text{[by (A.51) and (4.8)]} \\
&= \pi_1 (\eta + o(1)) \cdot c \log d^{-1} + O(d) \\
&\geq \pi_1 (\eta/2) \cdot c \log d^{-1} + o(c \log d^{-1}) > 0
\end{aligned}$$

for sufficiently small $c, d$, again a contradiction. Thus, (A.51) cannot hold either, showing that $N_1^* \sim I_1^{-1} \log d^{-1}$ and completing the proof. $\square$

**Acknowledgments.** The author would like to thank Gary Lorden for suggesting this problem and for his enthusiastic support of the author's thesis research, of which this work is an outgrowth. Amir Dembo and Tze Lai are also thanked for helpful discussions, as well as the Associate Editor and two referees, whose comments improved the paper.

DEPARTMENT OF MATHEMATICS
UNIVERSITY OF SOUTHERN CALIFORNIA
3620 SOUTH VERMONT AVENUE, KAP 108
LOS ANGELES, CALIFORNIA 90089
USA
E-MAIL: bartroff@usc.edu